\documentclass[11pt]{amsart}

\usepackage{amsmath}
\usepackage{amssymb}
\usepackage{graphicx}

\newtheorem{theorem}{Theorem}[section]
\newtheorem{corollary}[theorem]{Corollary}
\newtheorem{lemma}[theorem]{Lemma}

\theoremstyle{definition}
\newtheorem{definition}[theorem]{Definition}
\newtheorem{example}[theorem]{Example}
\newtheorem{remark}[theorem]{Remark}

\newcommand{\rav}{:\,\!=}
\newcommand{\ravref}[1]{\stackrel{(\ref{#1})}{=}}
\newcommand{\leqref}[1]{\stackrel{(\ref{#1})}{\leq}}
\newcommand{\geqref}[1]{\stackrel{(\ref{#1})}{\geq}}
\newcommand{\rref}[1]{$(\ref{#1})$}
\newcommand{\mm}[1]{{\mathbb{#1}}}
\numberwithin{equation}{section}

\title[On Lipschitz continuity]{On Lipschitz continuity of Value Functions\\ for Infinite Horizon Problem}

\author[D. Khlopin]{Dmitry Khlopin}

\address[D. Khlopin]{ Krasovskii Institute of Mathematics and Mechanics, Yekaterinburg, Russia;\ \  Ural Federal University, Yekaterinburg, Russia}
\email{{\tt khlopin@imm.uran.ru}}



\begin{document}

\begin{abstract}
    We investigate conditions of optimality for an infinite-horizon control problem and consider their correspondence with the value function.
 Assuming Lipschitz continuity of the value function, we prove that sensitivity relations plus  the normal form version of  the Pontryagin Maximum Principle is a necessary and sufficient condition for the optimality criteria that correspond to this value function. Different criteria of optimality under different asymptotic constraints may be used, including almost strong and classical optimality proposed by D.Bogusz.  Special attention is devoted to   weakly agreeable criteria.  We also obtain the conditions on control system (like controllability) that guarantee the Lipschitz continuity of the value function, without any other asymptotic conditions besides finiteness of the value function.

Some examples are discussed. In particular, it was shown that the same control, regarded as agreeable optimal and overtaking optimal control, can correspond to different (everywhere) value functions.


\end{abstract}

\maketitle

\section{Introduction}
Optimality conditions on control problems are usually constructed with the Pontryagin Maximum Principle (PMP)~\cite{ppp}. For infinite-horizon control problems, even in the free end-time case, the PMP can be degenerate  \cite{kr_as,Halkin}.

For infinite-horizon problems, one can often prove the nondegeneracy (the normality) of the PMP under additional assumptions. One may use uniform coercivity of running cost \cite[Corollary 4.1]{sagara}, strong convexity of the Hamiltonian \cite{tan_rugh},   or concavity of the dynamics function \cite[Hypothesis 5.1]{sagara}; or, one may impose asymptotic estimates on motions \cite{av,tauchnitz2015pontryagin} and  costate arcs \cite[Theorem~4]{kab},\cite[Remark~8]{kh_JDCS},\cite[Proposition~4]{Optim}, e.g.  their  total variation \cite{kr_as,au,bc,baumeister,sagara,norv}. The nondegeneracy of the PMP follows from the Lipschitz continuity of the value function  (see  \cite[Theorem~5.1]{kr_as}, \cite{JDCSnew,ye}) under the Michel condition (see \cite{michel}).
Besides, the Lipschitz continuity of the value function is also guaranteed by estimates on motions
and  costate arcs  e.g.  their  total variation, see \cite{kr_as,kab,av,bc,baumeister,sagara,ye}.

Through Dynamic Programming Principle (DPP), we prove that, for  a Lipschitz continuous value function,  the normal form of  PMP with sensitivity relations is a necessary and sufficient condition of optimality in view of this value function. To make it a straightforward consequence of the results \cite{fr2}  for finite horizon, we restate the optimality in  terms of the value function. Choosing a corresponding value function, we  obtain necessary conditions for  various optimality criteria under different asymptotic constraints, including the optimality in weighted Lebesgue spaces \cite{au,picken2015,tauchnitz2015pontryagin}, agreeable \cite{HaurieSethi,car1}, almost strong, and classical  \cite{slovak,picken} optimalities.
In these results, we only assume the corresponding value function to exist and to be Lipschitz continuous.
Generally, the same optimal control can correspond to several different value functions  (see  Example \ref{ex2}). For weakly agreeable  optimality criteria, we found a family of value functions parameterized by unboundedly increasing sequences. We also show  the conditions on a control system (see Theorems~\ref{11123},\ref{1112345}) that guarantee the Lipschitz continuity for a finite value function
   without an assumption  on the asymptotic behavior of trajectories or adjoint variables.

   We start with definitions, including the value function as a function satisfying DPP and optimal control corresponding to a given value function. In the next section,  we obtain a necessary (and sufficient) condition for such optimality. In Section~\ref{various}, we transfer  these results to different asymptotic constraints. Next, we study value functions for  agreeable optimality \cite{HaurieSethi,car1}. Section~\ref{appendix} is devoted to  the conditions on a control system that guarantee the Lipschitz continuity for a finite value function.
   The last section is devoted to examples.

 \section{Problem Statement and Definitions}

 {\bf The infinite horizon control problem.}
 Consider the following optimal control problem for infinite horizon:
\begin{eqnarray}
 \textrm{Minimize } \int_{0}^{+\infty} f_0(t,x,u)\, dt \label{sys0}\\
 \textrm{subject to } \dot{x}=f(t,x,u),\quad u\in P, \quad x\in \mm{R}^m; \label{sys}\\
 x(0)=b_*. \label{sysK}
\end{eqnarray}

Let $P$ be
a complete separable metric space.
 As for the class of admissible controls, we consider the set of all Lebesgue measurable
 functions $u:\mm{R}_{\geq 0}\to P$ that are bounded on every time compact and denote it by~${\mm{U}}$.

 Assume that $f:\mm{R}_{\geq 0}\times \mm{R}^m\times P\to\mm{R}^m$,
 $f_0:\mm{R}_{\geq 0}\times \mm{R}^m\times P\to\mm{R}$ are continuous and
 locally Lipschitz continuous in $x$;
  also, let
 $f$ satisfy the
 sublinear growth condition.
Then, to each $(b,t)\in \mm{R}^m\times \mm{R}_{\geq 0}$ and each  $u\in\mm{U}$,
we can assign a solution of \rref{sys} with the initial condition $x(t)=b$.
This solution is unique and it can be extended to the whole ${\mathbb{R}_{\geq 0}}$; denote it by~$x_{b,t,u}.$

 {\bf On the dynamic programming principle.}

For each $\theta\geq 0, y\in\mm{R}^m, T>\theta$, and $u\in\mm{U}$, set
$$J(\theta,y;u,T)\rav
 \int_{\theta}^T
 f_0\big(s,x_{y,\theta,u}(s),u(s)\big)\,ds.$$
Let
 $V^T(\theta,y)$ denote the infimum of the following problem in the interval $[\theta,T]$:
\begin{eqnarray*}
 \mm{P}(T)
  \left \{
 \begin{array}{rcl}
 \textrm{Minimize } \int_{\theta}^T  f_0(t,x,u)\, dt\\
 \textrm{subject to } \dot{x}=f(t,x,u),\quad u\in P, \quad x\in \mm{R}^m;\\
 x(\theta)=y.
 \end{array}            \right.
\end{eqnarray*}
\begin{definition}
For an interval $I\subset\mm{R}_{\geq 0}$
and a function ${V}: I\times \mm{R}^m \to \mm{R}\cup\{-\infty,+\infty\},$
we say that
${V}$
enjoys {\it DPP (Dynamic Programming Principle)} iff
$${V}(t,b)=\inf_{u\in\mm{U}} \Big[\int_{t}^\tau\!\! f_0(s,x_{b,t,u}(s),u(s))ds+{V}(\tau,x_{b,t,u}(\tau))\Big]\quad\forall b\in\mm{R}^m,
[t,\tau]\subset I.$$
\end{definition}
\begin{remark}\label{200}
For every $T>0$, $V^T: [0,T]\times \mm{R}^m \to \mm{R}$
satisfies DPP.
\end{remark}

 {\bf On conditions of optimality.}

Note that the improper integral in \rref{sys0} may not exist;  as a consequence,
 for control problems on infinite horizon,
  there are several optimality criteria 
  \cite{slovak,car1,cil,h_m,haur,HaurieSethi,picken,zz}. Hereinafter denote by $u^*$ the optimal control; however,
 we will always  specify which criterion do we use for $u^*$.
 Set also $x^*\equiv x_{b_*,0,u^*}$.
\begin{definition}
Let $\mm{V}:\mm{R}_{\geq 0}\times\mm{R}^m\to\mm{R}\cup\{-\infty,+\infty\}$ be a function enjoying
{DPP}.
We will say that
$u^*\in \mm{U}$
 {\it is optimal in view of} $\mm{V}$ iff, for all $T\geq 0$,
\begin{equation}\label{optimal}
 \mm{V}(T,x^*(T))+J(0,x^*(0);u^*,T)=\mm{V}(0,x^*(0)).
\end{equation}
\end{definition}
\section{Necessary and sufficient conditions  in terms of value functions.}

  Consider a Lipschitz continuous function~$h:\mm{R}^r\to\mm{R}$ and a point  $x\in\mm{R}^r$.
  The  Fr\'{e}chet  superdifferential  $\hat\partial^+ h(x)$  is  the  set  of  vectors $\zeta\in\mm{R}^r$ that satisfy
$$\limsup_{y\to x} \frac{h(y)-h(x)-\zeta(y-x)}{||y-x||}\leq 0.$$
  {\it The limiting superdifferential} $\partial^+ h(x)$ of $h$ at $x$ consists of all
$\zeta$ in $\mm{R}^r$  such that there exist sequences of $x_n\in{\mathbb{R}^r},\zeta_n\in \hat\partial^+ h (x_n)$  satisfying $x_n\to x,$ $h(x_n)\to h(x),$
$\zeta_n\to\zeta.$
For different (equivalent) definitions of the  limiting superdifferential, see \cite{borwein}.

{\bf The Pontryagin Maximum Principle.}

 Let the Hamilton--Pontryagin function
 $H$
 be given as follows:
 $$H(x,u,\psi,\lambda,t)\rav\psi f\big(t,x,u\big)-\lambda f_0\big(t,x,u\big).$$
 Let us introduce the relations
 \begin{eqnarray}
 \label{sys_psi}
 -\dot{\psi}(s)\!\!\!\!\!&=&\!\!\!\!\!\frac{\partial
 H}{\partial x}\big(x^*(s),u^*(s),\psi(s),\lambda,s\big),\\
 \label{maxH}
\sup_{v\in P}H\big(x^*(s),v,\psi(s),\lambda,s\big)\!\!\!\!\!&=&\!\!\!\!\!
 H\big(x^*(s),u^*(s),\psi(s),\lambda,s\big)\ a.e,\\
 -\psi(0)\!\!\!\!\!&\in&\!\!\!\!\! \partial^+_x \mm{V}(0,x^*(0)),\label{sens1}\\
  \big(\!H\big(x^*(s),u^*(s),\psi(s),1,s\big),
  \!-\psi(s)\big)\!\!\!\!\!&\in&\!\!\!\!\!
 \partial^+_{\ } \mm{V}(s,x^*(s))\ a.e.\label{sens2}
 \end{eqnarray}
 Here, $\partial^+_x \mm{V}$ is the limiting superdifferential of a map $\mm{R}^m\ni x\mapsto\mm{V}(s,x)\in\mm{R}.$

Halkin \cite{Halkin} proved that the Pontryagin Maximum Principle is a necessary condition of optimality for infinite-horizon problems:
 if an admissible $u^*\in\mm{U}$  is finitely optimal
\cite{Halkin} for problem \rref{sys0}--\rref{sysK}, then
 there exists a nontrivial solution $(\psi^*,\lambda^*)$ of system \rref{sys_psi}-\rref{maxH} for a.a. $s>0$; here, $\lambda^*\in\{0,1\}$.

 Remember that, for the corresponding control problem with finite horizon $T$,
under rather general conditions on the system, see, for example \cite{fr2,subbotina,vinter},
the Pontryagin Maximum Principle with $\lambda=1$ plus sensitivity relations \rref{sens1}--\rref{sens2} forms a necessary and sufficient condition of optimality.
Similar conditions arise in
          Hamilton-Jacobi-Isaacs PDEs for infinite horizon,
          see \cite{bc,baumeister}.

 Thereinafter in this section we also assume that $f$ and~$f_0$ are differentiable with respect to~$x$.

 {\bf Necessary conditions in terms of value function.}

 \begin{theorem}
\label{6767}
 Assume a locally Lipschitz continuous function $\mm{V}:\mm{R}_{\geq 0}\times\mm{R}^m\to\mm{R}$ satisfies DPP.
 Let  $u^*$ be optimal in view of   $\mm{V}$, i.e., satisfy \rref{optimal}.

Then, there exists a co-state arc $\psi\in C(\mm{R}_{\geq 0},\mm{R}^m)$ satisfying
PMP relations \rref{sys_psi}--\rref{maxH} with $\lambda= 1$ and
sensitivity relations \rref{sens1}--\rref{sens2}.
\end{theorem}
{\it Proof}
Fix  $T>0.$ Consider the following
 Bolza problem:
 \begin{eqnarray}
 \label{Bolza}
  \left \{
 \begin{array}{rcl}
 \textrm{Minimize } \int_{0}^T f_0(t,x,u)\, dt+\mm{V}(T;x(T))
 \\
 \textrm{subject to } \dot{x}=f(t,x,u),\quad u\in P, 
\\
x(0)=b. %
 \end{array} \right.
\end{eqnarray}
 Thanks to DPP, this problem has an optimal value, $\mm{V}(0,b).$
 Now, by \rref{optimal}, the control $u^*$ achieves the minimum of this problem, i.e.,
 $u^*$ is optimal for it.

 Thanks to \cite[Theorem~6.1]{fr2},
 for some co-state arc ${\psi}^T\in C([0,T],\mm{R}^m)$,
 PMP relations  \rref{sys_psi}--\rref{maxH} with $\lambda=1$
hold;
 moreover,  $\psi^T$
 satisfies \rref{sens2} (for a.a. $s\in[0,T]$) and \rref{sens1}.
 So, for every unbounded sequence of positive $T_n,$
 co-state arcs ${\psi}^{T_n}$ with $\lambda=1$
 satisfy \rref{sys_psi}--\rref{maxH} and  relations
 \rref{sens1},\rref{sens2} for a.a. $s\in[0,T_n].$

Remember that $\mm{V}$ is Lipschitz continuous, therefore
 $\partial^+_{x}\mm{V}(0,{x}^*(0))$ is a compact (see \cite{borwein}).
  Consider the sequence of ${\psi}^{T_n}(0)$.
 By \rref{sens1},
 this sequence has a limiting point $\zeta^*\in \partial^+_{x}\mm{V}(0,{x}^*(0))$.
 Consider a solution ${\psi}^*\in C(\mm{R}_{\geq 0},\mm{R}^m)$
  of \rref{sys_psi} with $\psi(0)=\zeta^*$; it satisfies \rref{sens1}.
 By the theorem on continuous dependence of a differential equation on its initial conditions,
 this solution is a partial limit of $\psi^{T_n}$ in the compact-open topology (on each time compact).
 Passing to the partial limit for a.a. positive $s$, one can provide that the co-state arc $\psi^*$ satisfies
 \rref{sys_psi}--\rref{maxH} and sensitivity relation \rref{sens2} for a.a. positive $s$.
 \qed

 {\bf Sufficient conditions  in terms of value function.}

 \begin{theorem}
\label{6766}
  Assume that a locally Lipschitz continuous function $\mm{V}:\mm{R}_{\geq 0}\times\mm{R}^m\to \mm{R}$ satisfies DPP.
Let some co-state arc $\psi\in C(\mm{R}_{\geq 0},\mm{R}^m)$ satisfy
 PMP relations \rref{sys_psi}--\rref{maxH} with $\lambda= 1$ and
sensitivity relation \rref{sens2}.

 Then, $u^*$ is optimal in view of   $\mm{V}$, i.e., satisfies \rref{optimal}.
\end{theorem}

{\it Proof}
Fix $T>0$. Consider  Bolza problem \rref{Bolza}.
 By DPP, this problem has an optimal value: $\mm{V}(0,b).$
 Now, \rref{optimal} holds for $T$ iff $u^*$ achieves the minimum of this problem, i.e., if
 $u^*$ is optimal for this  Bolza problem  with $b=b_*$.

 By \cite[Theorem~6.2]{fr2}, ${u}^*$ is optimal in this  Bolza problem  with $b=b_*$
     iff
 there exists a co-state arc ${\psi}^{T}\in C([0,T],\mm{R}^m)$ with $\lambda=1$ that satisfies \rref{sys_psi}--\rref{maxH}
 and  sensitivity relation \rref{sens2} for a.e. $s\in[0,T]$.
 Set $\psi^T\equiv\psi^*|_{[0,T]}.$
 Since the co-state $\psi^*$ with $\lambda=1$ satisfies \rref{sys_psi}--\rref{maxH},\rref{sens2}, we obtain \rref{optimal}  for all $T>0.$ \qed

Note that, in these  theorems, the local Lipschitz continuity of $\mm{V}$ may hold not in $\mm{R}_{\geq 0}\times\mm{R}^m$ but in a neighborhood of the graph of $x^*.$ In particular, one can require $\mm{V}$ to be locally Lipschitz continuous in a strongly invariant (for \rref{sys}) neighborhood of the graph of $x^*$.

 Of course, one would like to relax the condition of Lipschitz  continuity of the value function. On the other hand, Theorem \ref{11123} provides this based on the finiteness of the value function if the optimal-time function \cite{Vilev} for this control system is Lipschitz continuous.

%
%
%
\section{Value Functions under Asymptotic Constraints.}

\label{various}

 {\bf Asymptotic Constraints.}

 For controls  $u',u''\in\mm{U}$, for each $T\in\mm{R}_{\geq 0}$, the concatenation $u'\diamond_T u'\in\mm{U}$ is as follows:
 $(u'\diamond_\tau u'')(t)\rav
 u'(t)$ if $t<\tau,$ and
  $(u'\diamond_\tau u'')(t)\rav u''(t)$,
 if $t\geq\tau.$
\begin{definition}
We
say that a multi-valued map $\Omega_\diamond:\mm{R}^m\times \mm{R}_{\geq 0}\rightrightarrows \mm{U}$ {\it induces asymptotic constraints} iff,
for all $
(b,t)\in\mm{R}^m\times \mm{R}_{\geq 0}$,
\begin{eqnarray}\label{1123}
\Omega_\diamond(b,t)=\{u\diamond_T u_1\,:\,u\in\mm{U},u_1\in\Omega_\diamond(x_{b,t,u}(T),T)\},\quad\forall
{T}>t.
\end{eqnarray}
\end{definition}

It is easy to see that each of the constant multi-valued maps $\mm{U}$, $B(\mm{R}_{>0},P)$, $L_p(\mm{R}_{\geq 0},P)\cap\mm{U}$ (if $1\leq p<+\infty$)  induces asymptotic constraints.

Let us offer another example of asymptotic constraints.
Fix a nonempty set $M\subset\mm{R}^m.$
For all $(b,t)\in\mm{R}^m\times\mm{R}_{\geq 0}$,  denote  by $\Omega_M(b,t)$ the set of all $u\in\mm{U}$
   such that $\varrho(x_{b,t,u}(s),M)\to 0$ as $s\to+\infty.$
   Necessary conditions of optimality for problems under varying choice of $M$ were considered, for example, in \cite{ppp,pereira}; for  exit-time  control  problems,  see \cite{motta2014}.
   Let us check that  $\Omega_M$ does also induce asymptotic constraints.
 Indeed, for all $(t,b)\in\mm{R}_{\geq 0}\times\mm{R}^m$, $T>t$, $u,u_1\in\mm{U}$,
  both the statement $u\diamond_T u_1\in\Omega_M(x_{b,t,u}(t),t)$ and the statement $u_1\in\Omega_M(x_{b,t,u}(T),T)$
are equivalent to the fact that $\varrho(x_{b_1,T,u_1}(s),M)\to 0$ as $s\to+\infty$; here,
  $b_1\rav x_{b,t,u}(T)$. Thus,  $\Omega_M$ induces asymptotic constraints.

{\bf On conditions of optimality under asymptotic constraints}

Define
$V^{\diamond}:\mm{R}_{\geq 0}\times\mm{R}^m\to\mm{R}\cup\{-\infty,+\infty\}$
as follows:
 for all $(t,b)\in\mm{R}_{\geq 0}\times\mm{R}^m,$ set  $V^\diamond(t,b)=+\infty$ if  $\Omega_\diamond(t,b)=\varnothing$; otherwise,
 $$V^\diamond(t,b)\rav
 \inf_{u\in\Omega_\diamond(b,t)} \liminf_{T\to+\infty} J(t,b;u,T).$$

 \begin{theorem}
\label{6770d}
Assume that $\Omega_\diamond$  induces asymptotic constraints.
 Suppose that the function $V^{\diamond}$ is finite. Then,
 \begin{enumerate}
     \item $V^{\diamond}$  enjoys the Dynamic Programming Principle (DPP);
     \item
 $u^*$ is optimal in view of $V^{\diamond}$ if $u^*$ lies in $\Omega_{\diamond}(b_*,0)$ and satisfies
\begin{equation}\label{464}
 \liminf_{T\to+\infty} {J}(0,b_*;u^*,T)=V^\diamond(0,b_*).
 \end{equation}
   \end{enumerate}
\end{theorem}
{\it Proof}
Note that, since  $V^\diamond$ is everywhere finite, $\Omega_{\diamond}$ is everywhere nonempty.
By the definition of $V^\diamond$,
for all $\tau>0,t\in[0,\tau[,b\in\mm{R}^m,$ we have
\begin{eqnarray*}
\begin{split}
V^{\diamond}(t,b)=&\inf_{u\in\Omega_\diamond(b,t)} \liminf_{T\to+\infty} \Big[J(t,b;u,\tau)+J(\tau,x_{b,t,u}(\tau);u,T)\Big]\\
\ravref{1123}&
\inf_{u'\in\mm{U},u''\in\Omega_\diamond(x_{b,t,u'}(\tau),\tau)} \Big[J(t,b;u'\diamond u'',\tau)\\
\ &+\liminf_{T\to+\infty} J(\tau,x_{b,t,u'\diamond u''}(\tau);u'\diamond u'',T)\Big]\\
 =&\inf_{u'\in\mm{U}} \Big[J(t,b;u',\tau)
\\ \ &+\inf_{u''\in\Omega_\diamond(x_{b,t,u'}(\tau),\tau)} \liminf_{T\to+\infty} J(\tau,x_{b,t,u'}(\tau);u'',T)\Big]\\
=&\inf_{u'\in\mm{U}} \Big[J(t,b;u',\tau)+V^{\diamond}(\tau,x_{b,t,u'}(\tau))\Big].
\end{split}
\end{eqnarray*}

Assume \rref{464} holds for $u^*\in\Omega_\diamond(b_*,0)$. Then, for some   unbounded sequence of positive $\tau_n$, we have
$\displaystyle\lim_{n\to\infty} {J}(0,b_*;u^*,\tau_n)=V^\diamond(0,b_*).$
Now,
\begin{eqnarray*}
V^\diamond(0,b_*)&\ravref{464}&
\lim_{n\to\infty}\big[J(0,b_*;u^*,T)+{J}(T,x^*(T);u^*,\tau_n)\big]\\
&=&J(0,b_*;u^*,T)+\lim_{n\to\infty}{J}(T,x^*(T);u^*,\tau_n)\\
&\geqref{1123}&J(0,b_*;u^*,T)+
 \inf_{u\in\Omega_\diamond(x^*(T),T)} \liminf_{n\to\infty} J(T,x_{b,t,u}(T);u,\tau_n)\\
&\geq& J(0,b_*;u^*,T)+{V}^\diamond(T,x^*(T))\qquad \forall T\geq 0.
\end{eqnarray*}
 Thanks to the Dynamic Programming Principle, it implies \rref{optimal} for $\mm{V}={V}^\diamond$.\qed

{\bf On conditions without asymptotic constraints.}

Define the function  $V^{inf}$ from $\mm{R}_{\geq 0}\times\mm{R}^m$ to $\mm{R}\cup\{-\infty,+\infty\}$
as follows:
\begin{eqnarray*}
 V^{inf}(t,b)=
 \inf_{u\in\mm{U}} \liminf_{T\to\infty} J(t,b;u,T)\qquad \forall (t,b)\in\mm{R}_{\geq 0}\times\mm{R}^m.
\end{eqnarray*}
One could easily prove that $V^{inf}\equiv V^\diamond$ if $\Omega_\diamond\rav \mm{U}$.
Now, Theorems~\ref{6767} and~\ref{6770d} imply the corresponding necessary conditions for $u^*$.
\begin{corollary}
\label{6770}
 Assume that $f$ and~$f_0$ are differentiable with respect to $x$.

 Let $V^{inf}$ be finite and   locally Lipschitz continuous.
 Let $u^*\in\mm{U}$  satisfy
 $$\displaystyle \liminf_{T\to+\infty} {J}(0,b_*;u^*,T)=V^{inf}(0,b_*).$$

Then,
 there exists a co-state arc $\psi\in C(\mm{R}_{\geq 0},\mm{R}^m)$ satisfying  PMP relations \rref{sys_psi}--\rref{maxH} with $\lambda= 1$ and
 sensitivity relations \rref{sens1}--\rref{sens2} for $\mm{V}=V^{inf}$.
\end{corollary}
 Note that $u^*$ is  optimal in view of $V^{inf}$ if $u^*$ is an overtaking optimal control \cite{HaurieSethi,car1}. Therefore, we obtain necessary conditions of overtaking optimality.
 These necessary conditions, including the effective transversality condition at infinity, can be found in \cite{av,belyakov}.
   However, the necessary conditions in \cite{av,belyakov} exploit the assumptions that are not required for Corollary~\rref{6770}; however, Corollary~\ref{6770} assumes the value function to be known.

Analogously, set $\Omega_\diamond\equiv L_p(\mm{R}_{\geq 0},P)\cap\mm{U}$. If $u^*$ achieves $V^{\diamond}(0,b_*)$, \rref{sens1} holds. For the linear case, this sensitivity relation  was proved in \cite[(7.13)]{au}.

{\bf On almost strong and classical optimalities.}

Recall the optimality criteria from \cite{slovak} and \cite{picken}.
  For each $(b,t)\in\mm{R}^m\times\mm{R}$,  denote
    by $\Omega_{\mathcal{R}}(b,t)$ ($\Omega_{\mathcal{L}}(b,t)$) the set of all $u\in\mm{U}$
   such that a map $[t,+\infty[\ni s\mapsto f_0(s,x_{b,t,u}(s),u(s))$  has an improper Riemann (Lebesgue) integral.
We claim that $\Omega_{\mathcal{R}}$  induces asymptotic constraints.
 Fix  $(b,t)\in\mm{R}^m\times\mm{R}_{\geq 0}$, $u,u_1\in\mm{U}$, and $T>t.$ Set $b_1\rav x_{b,t,u}(T)$.
By definition, $u_1\in\Omega_{\mathcal{R}}(b_1,T)$ iff, for all $T>0$, there exists
 a  limit of  $\int_{T}^\tau f_0(s,x_{b_1,T,u_1}(s),u_1(s))\,ds$,
   i.e., if, specifically,
   $\int_{t}^\tau f_0(s,x_{b_1,T,u\diamond_T u_1}(s),(u\diamond_T u_1)(s))\,ds$ has a limit as $\tau\to+\infty$.
 By $x_{b_1,T,u\diamond_T u_1}=x_{b,t,u\diamond_T u_1}$,
 we obtain ``$\supseteq$'' in \rref{1123} for $\Omega_{\diamond}\equiv\Omega_{\mathcal{R}}$.
 On the other hand, $u\in\Omega_{\mathcal{R}}(b_1,T)$ if $u\in\Omega_{\mathcal{R}}(b,t)$; thanks to
  $u=u\diamond_T u,$ we also obtain the converse inclusion.
The proof for $\Omega_{\mathcal{L}}$  is similar.

Let the  mappings $V^{\mathcal{L}}$,
$V^{\mathcal{R}}$ from $\mm{R}_{\geq 0}\times\mm{R}^m$ to $\mm{R}\cup\{-\infty,+\infty\}$
be $V^{\diamond}$ with $\Omega_{\diamond}=\Omega_{\mathcal{L}}$, $\Omega_{\diamond}=\Omega_{\mathcal{R}}$, respectively.
Under the assumption $V^{\mathcal{L}}(0,x^*(0))\in\mm{R},$
it is easy to verify that $u^*$ is {\it classical optimal}
\cite[Definition 7.5]{slovak}, \cite[(L1)]{picken} iff $u^*\in\Omega_{\mathcal{L}}(x^*(0),0)$ holds \rref{464} for
$V^{\diamond}=V^{\mathcal{L}}.$
Similarly,
under the assumption  $V^{\mathcal{R}}(0,b_*)\in\mm{R}$,
$u^*$ {\it is almost strongly optimal} \cite[Definition 7.8]{slovak},\cite[(R1)]{picken} iff
$u^*\in\Omega_{\mathcal{R}}(x^*(0),0)$ holds \rref{464} for $V^{\diamond}=V^{\mathcal{R}}$.
Theorems~\ref{6767} and~\ref{6770d} imply necessary  conditions of optimality for these criteria.
In \cite[Sect. 7]{slovak}, see rather general conditions of existence for such optimal controls.

\section{On Agreeable Control}

\begin{definition}
Call $u^*\in\mm{U}$
{\it a weakly agreeable} control  \cite[Definition 3.2(iii)]{car1}
iff   this control (with its motion $x^*$, $x^*(0)=b_*$)  satisfies
\[ \liminf_{T\to+\infty}
 \Big(
 J\big(0,b_*;u^*,t)+V^T\big(t,x^*(t)\big)
 -
 \inf_{u\in{\mm{U}}}\,J\big(0,b_*;u,T\big)
 \Big)
 \leq 0,\quad \forall t\geq0.\]
\end{definition}

\begin{definition}
Call $u^*\in\mm{U}$
 {\it an agreeable} control \cite[Definition 3.2(ii)]{car1} iff
     this control (with its motion $x^*$, $x^*(0)=b_*$) satisfies
\[ \lim_{T\to+\infty}
 \Big(
 J(0,b_*;u^*,t)+V^T(t,x^*(t))
 -
 \inf_{u\in{\mm{U}}} J(0,b_*;u,T)
 \Big)
 \leq 0,\quad \forall t\geq0.\]
\end{definition}

\begin{lemma}
\label{777}
  For a
  control $u^*\in\mm{U}$ with its motion $x^*$, $u^*$ is weakly agreeable iff
  there exists an unbounded sequence of positive
  $\tau_n$ such that, for all $T>0,$
\begin{eqnarray}
 {J}(0,x^*(0);u^*,T)=\lim_{n\to\infty} \big[V^{\tau_n}(0,x^*(0))-V^{\tau_n}(T,x^*(T))\big]
 \label{op_uow}
\end{eqnarray}
\end{lemma}
{\it Proof}
 Let $u^*\in\mm{U}$ be weakly agreeable.
 Then, for a natural $n$, we can choose $\tau_n>n$ such that
${J}(0,x^*(0);u^*,n)+V^{\tau_n}(n,x^*(n))- V^{\tau_n}(0,x^*(0))<1/n.$
 For all  $T\in[0,n],$ the relation $u^*\diamond_{n} \mm{U}\subset u^*\diamond_{T} \mm{U}$ holds. Hence,
\begin{eqnarray*}
 {J}(0,x^*(0);u^*,T)+ V^{\tau_n}(T,x^*(T))\leq{J}(0,x^*(0);u^*,n)+ V^{\tau_n}(n,x^*(n)),
\end{eqnarray*}
 also, by DPP,
 $0\leq {J}(0,x^*(0);u^*,T)+V^{\tau_n}(T,x^*(T))- V^{\tau_n}(0,x^*(0))<1/n$.
 Passing to the limit as $n\to\infty$, we obtain \rref{op_uow} for this sequence of $\tau_n.$

 The converse implication is clear.\qed

{\bf Necessary conditions of weakly agreeable optimality}

 \begin{corollary}
\label{6769}
 Let $u^*\in\mm{U}$ be a
 weakly agreeable control
 in \rref{sys0}--\rref{sysK}, i.e.,
 for some unbounded sequence $\tau_n\uparrow\infty$, let $u^*$ satisfy \rref{op_uow} for all $T>0$.

Let a function $V^{\infty}:\mm{R}_{\geq 0}\times \mm{R}^m\to \mm{R}\cup\{-\infty,+\infty\}$, defined as
\begin{eqnarray}
\label{tauinfty}
 V^\infty(t,b)\rav\liminf_{n\to\infty} V^{\tau_n}(t,b),\quad \forall (t,b)\in\mm{R}_{\geq 0}\times\mm{R}^m,
\end{eqnarray}
be finite and locally Lipschitz continuous.

Then, there exists a co-state arc $\psi\in C(\mm{R}_{\geq 0},\mm{R}^m)$ such that  PMP relations \rref{sys_psi}--\rref{maxH} for $\lambda= 1$ and
 sensitivity relations \rref{sens1}--\rref{sens2} for $\mm{V}=V^{\infty}$ hold.
\end{corollary}
{\it Proof}
 Fix  $T>0,(t,b)\in [0,T]\times\mm{R}^m.$
 By Remark~\ref{200}, $V^{\tau_n}$ satisfies DPP.
 Passing to the lower limit as $n\to\infty$, we have
 $$V^\infty(t,b)=\inf_{u\in\mm{U}} \Big[\int_{t}^T f_0(s,x_{b,t,u}(s),u(s))ds+V^{\infty}(T,x_{b,t,u}(T))\Big]$$
 for all $T>0,(t,b)\in [0,T]\times\mm{R}^m.$
 Thus, $V^\infty$ satisfies DPP.

Fix  $T>0$. By \rref{op_uow}, there exists a converging to $0$ sequence of  $\varepsilon_n$ such that
$${J}(0,x^*(0);u^*,T)+V^{\tau_n}(T,x^*(T))= V^{\tau_n}(0,x^*(0))+\varepsilon_n$$
holds.
Passing to the lower limit as $\tau_n\to\infty$, we obtain
$${J}(0,x^*(0);u^*,T)+V^{\infty}(T,x^*(T))= V^{\infty}(0,x^*(0)),\qquad \forall T>0.$$
 Now, Theorem~\ref{6767} with $\mm{V}= V^{\infty}$ completes the proof.\qed

{\bf Sufficient conditions of weakly agreeable optimality}

\begin{lemma}
\label{560}
 Let $\tau$ be an unbounded sequence of positive numbers.
 Let the value $V^\infty(0,b_*)$, defined in \rref{tauinfty}, be finite.
 Let  $u^*$ satisfy \rref{optimal} for  $\mm{V}=V^\infty$.

 Then,  $u^*$ is weakly agreeable, and \rref{op_uow} holds for some subsequence of $\tau$.
 \end{lemma}
{\it Proof}
 Note that, by \rref{optimal}, the finiteness of $V^\infty(0,x^*(0))$ implies that $V^\infty(t,x^*(t))$ is the same for all $t>0$.
 Now,  there exists a subsequence of $\tau'_k\rav\tau_{n(k)}$
 such that
 $V^{\tau'_k}(k,x^*(k))\leq V^\infty(k,x^*(k))+1/k$ for a natural $k$.
Fix these $\tau'_{k}.$

 For each  $k\in\mm{N},T\in[0,k],$  we have $u^*\diamond_{k} \mm{U}\subset u^*\diamond_{T} \mm{U}$. Hence,
\begin{eqnarray*}
 {J}(0,b_*;u^*,T)+ V^{\tau'_k}(T,x^*(T))&\leq&{J}(0,b_*;u^*,k)+ V^{\tau'_k}(k,x^*(k))\\
  &\leq& {J}(0,b_*;u^*,k)+V^\infty(k,x^*(k))+1/k\\
  &\ravref{optimal}&
  V^\infty(0,x^*(0))+1/k\\
  &=& V^\infty(T,x^*(T))+{J}(0,b_*;u^*,T)+1/k.
\end{eqnarray*}
 Thus,
 $V^{\tau'_k}(T,x^*(T))$ converges to $V^\infty(T,x^*(T))$ for all $T\geq 0$ as $k\to\infty.$

 Now, for each  $T>0$, passing to the limit, we obtain
\begin{eqnarray*}
  0&\ravref{optimal}&J(0,x^*(0);u^*,T)+{V}^{\infty}(T,x^*(T))-{V}^{\infty}(0,x^*(0))\\
  &=&
  \lim_{k\to\infty}\big[J(0,x^*(0);u^*,T)+{V}^{\tau'_k}(T,x^*(T))-{V}^{\tau'_k}(0,x^*(0))\big].
\end{eqnarray*}
  Thus,  \rref{op_uow} holds for  $\tau'$, and, by Lemma~\ref{777}, $u^*$ is weakly agreeable.
 \qed

Applying Lemma~\ref{560} and Theorem~\ref{6766}, we obtain a sufficient condition for  a weakly agreeable control.
\begin{corollary}
\label{561}

For some unbounded sequence $\tau\uparrow\infty,$ let $V^\infty$ be finite and locally Lipschitz continuous.
Let a co-state arc $\psi\in C(\mm{R}_{\geq 0},\mm{R}^m)$ satisfy
 PMP relations \rref{sys_psi}--\rref{maxH} with $\lambda= 1$ and
sensitivity relation \rref{sens2} for $\mm{V}=V^\infty$.

 Then,  \rref{op_uow} holds for some subsequence of $\tau$, and $u^*$ is weakly agreeable.
 \end{corollary}

See the sufficient conditions  of agreeable optimality
 in terms of the asymptotic behavior of $\psi(s)(x(s)-x^*(s))$ for large~$s$ in \cite[Theorems 2.5 and 2.6]{haur}.
 Such conditions are often used in proofs of the turnpike property \cite{zz}.

{\bf Conditions for agreeable optimality.}

 \begin{theorem}
\label{676767_}
Let a function ${V}^{all}$, defined as follows,
\begin{eqnarray}\label{4200}
{V}^{all}(t,b)\rav\lim_{T\to+\infty} V^{T}(t,b), \qquad \forall (t,b)\in\mm{R}_{\geq 0}\times\mm{R}^m,
\end{eqnarray}
be well-defined and locally Lipschitz continuous.

 Then,
 the following conditions are equivalent:
 \begin{enumerate}
   \item $u^*$ is weakly agreeable;
  \item $u^*$ is agreeable;
  \item
  there exists a co-state arc $\psi\in C(\mm{R}_{\geq 0},\mm{R}^m)$ satisfying PMP relations \rref{sys_psi}--\rref{maxH} with $\lambda= 1$ and
 sensitivity relations \rref{sens1},\rref{sens2} for $\mm{V}=V^{all}$.
 \end{enumerate}
\end{theorem}
{\it Proof}  $2)\Rightarrow 1)$ was proved in \cite[Proposition~3.2]{car1}.

$1)\Rightarrow 3).$
By Lemma~\ref{777},  $u^*$ satisfies
\rref{op_uow}  for some unbounded sequence of positive $\tau_n.$
Note that $V^\infty\equiv V^{all},$ and all conditions of Corollary~\ref{6769} hold.
Therefore,  for some co-state arc $\psi\in C(\mm{R}_{\geq 0},\mm{R}^m)$, relations \rref{sys_psi}--\rref{maxH} for $\lambda= 1$ and
  relations \rref{sens1}--\rref{sens2} for $\mm{V}=V^{\infty}$ hold.
By $V^\infty\equiv V^{all},$  $1)\Rightarrow 3)$ is proved.

 $3)\Rightarrow 2).$
Assume the contrary.
Then, by the definition of  agreeable controls,
there exist $t\geq 0$ and
an unbounded sequence of positive $\tau_n$ such that
\[ \liminf_{n\to\infty}
 \Big(
 J(0,b_*;u^*,t)+V^{\tau_n}(t,x^*(t))
 -
 V^{\tau_n}(0,x^*(0))
 \Big)
 >0.\]
 On the other hand,
 $V^\infty\equiv V^{all}$ and all conditions of Corollary~\ref{561} hold, and
 \rref{op_uow} holds for $u^*$ with every unbounded sequence of positive numbers.
 However, it contradicts the choice of $\tau_n.$\qed

 Some conditions of validity of $V^{all}$ were shown in  \cite[Theorem 3.3]{dalio}. In this case,  $V^{\inf}\equiv V^{all}$ holds  \cite[(3.13)]{dalio}.
 Then, a weakly agreeable control is also  weakly overtaking optimal \cite{car1}.
 Conditions \cite[Hypothesis 3.1(i)-(iv), Hypothesis A.1]{sagara} guarantee the Lipschitz continuity for well-defined $V^\infty$ and $V^{all}$
 (it is sufficient to repeat the proof of \cite[Theorem~A.1]{sagara} verbatim).
 On the other hand, in Example~\ref{ex2}, the value functions   $V^{\inf}$ and~$V^{all}$ are Lipschitz continuous and there exists a control that is optimal for both value functions, however, $V^{\inf}>V^{all}$ holds everywhere.

\section{On Conditions of Lipschitz Continuity of Value Functions}
 Usually, the Lipschitz continuity of a value function is guaranteed by asymptotic conditions on $f,f_0,J$ or on solutions of \rref{sys},\rref{sys_psi}, see \cite{kr_as,kab,av,bc,baumeister,sagara,ye}.
 We will use the Lipschitz continuity \cite{Vilev,MottaRampazzo} of the optimal-time function  of the control system.
 Let us obtain the conditions to guarantee the Lipschitz continuity of the value function in absence of any asymptotic conditions besides the finiteness of this function.

    Consider nonnegative integers $r,s$ ($r+s=m$), a set $\mm{W}\subset\mm{R}^r$, and  a compact subset $P'\subset P$.  Consider optimal-time function $Q_{\mm{W}}$,
    optimal-time function $Q'_{\mm{W}}$ under additional condition $u(t)\in P'(\forall t\geq 0),$ and
    maximal-time function $Q^{\mm{W}}$ as follows: for all $y'=(w',z')\in{G},z\in \mm{R}^s,t'\geq 0$,
\begin{equation*}
\begin{split}
Q'_{\mm{W}}(t',y',z)\rav& \inf\!\big\{\tau\geq 0\,:\,\!\exists\ u\in\mm{U}\cap B(\mm{R}_{\geq 0},P'),  x_{y',t',u}(t'\!+\!\tau)\in \mm{W}\!\times\!\{z\}\big\},\\
Q_{\mm{W}}(t',y',z)\rav& \inf\!\big\{\tau\geq 0\,:\,\!\exists\ u\in\mm{U},  x_{y',t',u}(t'\!+\!\tau)\in \mm{W}\!\times\!\{z\}\big\},\\
Q^{\mm{W}}(t',y',z)\rav& \sup\!\big\{\tau\geq 0\,:\,\!\exists\ u\in\mm{U},  x_{y',t',u}(t'\!+\!\tau)\in \mm{W}\!\times\!\{z\}\big\},
\end{split}
\end{equation*}
 $Q'_{\mm{W}}(y',z,t')\rav+\infty$ if no motion joins $(t',y')$ to $[t',\infty[\times \mm{W}\times\{z\}$  under the requirement $u(t)\in P'$ for all $t\geq t'$,
 $Q'_{\mm{W}}(y',z,t')\rav+\infty$ if no motion joins $(t',y')$ to $[t',\infty[\times \mm{W}\times\{z\}$, and
  $Q^{\mm{W}}(t',y',z)\rav +\infty$ if there exists a motion that does not join  $(t',y')$ to $[t',\infty[\times \mm{W}\times\{z\}$.

Let $I,W$,  and $Z$ be non-empty open subsets of $\mm{R}_{\geq 0},int\,\mm{W},\mm{R}^s$ respectively.

We will prove the following theorems:


 \begin{theorem}
\label{11123}
  Let $\mm{V}:\mm{R}_{\geq 0}\times\mm{R}^m\to\mm{R}$ enjoy DPP.
   Assume that  a function $S:Z\to \mm{R}$ and a positive locally Lipschitz continuous function $R:I\times \mm{W}\to \mm{R}_{>0}$ satisfy
  $$\mm{V}(t,(w,z))=R(t,w) S(z) \qquad \forall (t,w,z)\in I\times \mm{W}\times Z.$$

   For a point of $I\times W\times Z$, let there exist a neighborhood $I'\times W'\times Z'$ of this point and a positive $L$ such that
   $$Q'_{\mm{W}}(t',(w',z'),z)\leq L||z'-z||\quad \forall t'\in I',w'\in W',z,z'\in Z'.$$

Then, $\mm{V}$ is
 locally Lipschitz continuous on $I\times W\times Z$.
\end{theorem}

 \begin{theorem}
\label{1112345}
  Let $\mm{V}:\mm{R}_{\geq 0}\times\mm{R}^m\to\mm{R}$ enjoy DPP.
   Assume that  a function $S:Z\to \mm{R}$ and a positive locally Lipschitz continuous function $R:I\times \mm{W}\to \mm{R}_{>0}$ satisfy
  $$\mm{V}(t,(w,z))=R(t,w) S(z)\qquad \forall  (t,w,z)\in I\times \mm{W}\times Z.$$

   For a point of $I\times W\times Z$, let there exist a neighborhood $I'\times W'\times Z'$ of this point and a positive $L$ such that
   $$||z'-z||/L\leq Q_{\mm{W}}(t',(w',z'),z)\leq Q^{\mm{W}}(t',(w',z'),z)\leq L||z'-z||$$ holds for all $t'\in I',w'\in W',z,z'\in Z'$.

Then, $\mm{V}$ is  locally Lipschitz continuous on $I\times W\times Z$.
\end{theorem}
 \begin{corollary}
\label{111234}
  Assume  that $\mm{V}:\mm{R}_{\geq 0}\times\mm{R}^m\to\mm{R}$ satisfies DPP.

   Assume that
   $$0\in int\,conv\,\{f(t,x,u)\,:\,u\in P'\}\quad \forall (t,x)\in G$$
    holds for some open set $G\subset\mm{R}_{\geq 0}\times\mm{R}^m.$

   Assume that a function $S:\mm{R}^m\to\mm{R}$ and a locally Lipschitz continuous function $R:G\to \mm{R}_{>0}$ satisfy    $$\mm{V}(t,x)=R(t) S(x) \qquad\forall  (t,x)\in G.$$

Then, $\mm{V}$ is
 locally Lipschitz continuous in $G$.
\end{corollary}
 \begin{corollary}
\label{11123456}
  Assume  that $\mm{V}:\mm{R}_{\geq 0}\times\mm{R}^m\to\mm{R}$ satisfies DPP.

   Assume that
    $$0\not\in cl\,conv\,\{f(t,x,u)\,:\,(t,x)\in G,u\in P\}$$ 
   !holds!
   for some open set $G\subset\mm{R}_{\geq 0}\times\mm{R}^m$.

   Assume that a function $S:\mm{R}^m\to\mm{R}$ and a locally Lipschitz continuous function $R:G\to \mm{R}_{>0}$ satisfy  
    $$\mm{V}(t,x)=R(t) S(x) \qquad\forall  (t,x)\in G.$$

Then, $\mm{V}$ is
 locally Lipschitz continuous in $G$.
\end{corollary}


{\bf Proof\ of Theorem~\ref{11123}.}

 Fix $(t_*,b_W,b_Z)\in I\times W\times Z.$
 Take a fitting neighborhood $I'\times W'\times Z'$ of this point with some $L>1$. It is safe to assume $cl\,(I'\times W'\times Z')\subset I\times W\times Z$ to be compact, !$I'$ be an interval!.
 It will suffice to prove that
  $\mm{V}$ is Lipschitz continuous in another, possibly smaller neighborhood of $(t_*,b_W,b_Z)$.
   Set $b\rav(b_W,b_Z)$, $G'\rav W'\times Z'.$

     Thanks to the sublinear growth of $f,$ there exist a compact $K$ ($cl\,G'\subset int\,K\subset \mm{W}\times Z$) and $\gamma>0$ such that $y\in G'$
  implies
  $$x_{y,\tau,u}(t+\tau)\in K,t+\tau\in I\quad \forall
  u\in\mm{U}\cap B(\mm{R}_{\geq 0},P'), t\in I',\tau\in[0,\gamma].$$
    Now, we can choose
  $M>\max\{L+1,|\mm{V}(t_*,b)|\}$ such that
  \begin{eqnarray*}
  \begin{split}
  f_0(t+\tau,y,u)<M,\quad  &|R(t+\tau,\bar{w})|\in[1/M,M],\\
  ||f(t+\tau,y,u)||<M,\quad
  &|R(t,w')-R(t+\tau,\bar{w})|\leq M||w'-w||+M\tau
  \end{split}
  \end{eqnarray*}
   hold for all $y\in K,(w,z),(w',z')\in K\cap(\mm{W}\times Z),u\in P',t\in I,\tau\in[0,\gamma]$. Now, we obtain $||x_{y,t,u}(t+\tau)-y||\leq M\tau$ for all
  $y\in G',\tau\in I.$ Also, we have $|S(b_Z)|\leq |\mm{V}(t_*,b)|/R(t_*,b_W)\leq M^2$.
    Decreasing the neighborhood $G'\rav W'\times Z'$ of $b$, we can assume $diam\, Z'<\gamma/M.$

 Fix $t'\in I.$  Consider $w'\in W',z',z''\in Z'$.  Set $y'\rav (w',z')\in G'.$
 By the choice of  $L$, we have
 $$Q'_{\mm{W}}((w',z'),z'',\tau)\leq L||z'-z''||<M\,diam\, Z'<\gamma.$$
 By the definition of $Q_{\mm{W}}$, we can choose
    $u\in\mm{U}\cap B(\mm{R}_{\geq 0},P')$, $\bar{w}\in \mm{W}$, and $\tau>0$ enjoying $$x_{y',t',u}(t'+\tau)=(\bar{w},z''),\quad
   \tau \leq (L+1)||z'-z''||\leq M\,diam\, Z'<\gamma.$$
   By $\tau<\gamma$ and  $y'\in G'$, we
   obtain  $(\bar{w},z'')\in K$ and
 \begin{eqnarray}\label{1047}
   ||\bar{w}-w'||\leq ||x_{y',t',u}(t'+\tau)-y'||\leq M\tau< M^2||z'-z''||.
\end{eqnarray}
  Set $t''\rav t'+\tau.$
   By  DPP, we have $$\mm{V}(t',y')\leq J(t',y';u,t'')+\mm{V}(t'',x_{y',0,u}(t'')).$$
   In view of $\mm{V}(t,(w,z))=R(t,w) S(z)$ and
    $f_0(t,y,u)\leq M,$ we obtain
 \begin{equation}\label{1048}
 \begin{split}
   R(t',w') S(z')&=\mm{V}(t',y')\\
   &\leq
   M\tau+
   R(t'',\bar{w}) S(z'')=M\tau+\mm{V}(t'',(\bar{w},z'')).
  \end{split}
\end{equation}

   In particular, 
   in the case  $y'\rav y,z''\rav b_Z$, we can provide $(\bar{w},z'')\in K$, $\tau<\gamma$, and
  $$R(t',w') S(z')\leqref{1048}
   R(t'',\bar{w}) S(b_Z)+M\tau\leq M |S(b_Z)|+M\leq 2M^3$$
   by the choice of $M$.
   In the case  $y'\rav b,z''\rav z$, we have $(\bar{w},z)\in K$, $\tau<\gamma$ and, at last,
   $$R(t'',w'')S(z)\geqref{1048} R(t',b_W) S(b_Z)-\tau M \geq -M^3-\tau M\geq -2M^3.$$
   According  to $|R|>1/M$, we also have $|S(z)|\leq 2M^4$ for all $(w,z)\in G.$

  Return to general case.   Consider every  $w'\in W',z',z''\in Z'$  and, additionally, every $w''\in W',t\in I$. Then, we again obtain  $(\bar{w},z'')\in K$ and $\tau<\gamma$ satisfying \rref{1047} and \rref{1048}. Now, by \rref{1048},
  for all $t,t'\in I$,$w',w''\in W',$ and $z',z''\in Z',$
   $\mm{V}(t',w',z')-\mm{V}(t,w'',z'')$ does not exceed
\begin{eqnarray*}
   M\tau+               (R(t'\!+\!\tau,\bar{w})-    R(t',w')+R(t',w')-R(t,w''))S(z'')\\
   \leq M\tau+ 2M^5(\tau+||\bar{w}-w'||+||{w}''-w'||+|t'-t|)\\
   \leqref{1047} (M+ 2M^5+2M^6)\tau+2M^5(||{w}''-w'||+|t'-t|)\\
   \leqref{1047} 
 7M^7(||z'-z''||+||w''-w'||+|t'-t|).
\end{eqnarray*}
 Thus,  $$|\mm{V}(t',w',z')-\mm{V}(t,w'',z'')|\leq 7M^7(||z'-z''||+||w''-w'||+|t'-t|)$$ holds
 for all $(t',w',z'),(t,w'',z'')$ from
  some neighborhood $I\times W'\times Z'$ of each $(t_*,b_W,b_Z)\in I\times W\times Z$.
  \qed

\label{appendix}
{\bf Proof\ of Theorem~\ref{1112345}.\ }

 Fix $(t_*,b_W,b_Z)\in I\times W\times Z.$
 Take a fitting neighborhood $I'\times W'\times Z'$ of this point with some $L>1$. It is safe to assume $cl\,(I'\times W'\times Z')\subset I\times W\times Z$ to be compact,  !$I'$ be an interval!.
 It will suffice to prove that
  $\mm{V}$ is Lipschitz continuous in another, possibly smaller neighborhood of $(t_*,b_W,b_Z)$.
   Set $b\rav(b_W,b_Z)$, $G'\rav W'\times Z'.$

     Thanks to the sublinear growth of $f,$ there exist a compact $K$ ($cl\,G'\subset int\,K\subset \mm{W}\times Z$) and $\gamma>0$ such that 
  $$x_{y,\tau,u}(t+\tau)\in K,t+\tau\in I\qquad \forall y\in G',\tau\in[0,\gamma],
 u\in\mm{U}\cap B(\mm{R}_{\geq 0},P'),t\in I'.$$
    Now, we can choose
  $M>\max\{L+1,|\mm{V}(t_*,b)|\}$ such that
  \begin{eqnarray*}
  \begin{split}
  |f_0(t+\tau,y,u)|<M,\ &|R(t+\tau,\bar{w})|\in[1/M,M],\\
  ||f(t+\tau,y,u)||<M,\ &|R(t,w')-R(t+\tau,\bar{w})|\leq M||w'-w||+M\tau
  \end{split}
  \end{eqnarray*}
   hold for all $y\in K,(w,z),(w',z')\in K\cap(\mm{W}\times Z),u\in P',t\in I,\tau\in[0,\gamma]$. Now, we obtain $||x_{y,t,u}(t+\tau)-y||\leq M\tau$ for all
  $y\in G',\tau\in I.$ Also, we have $|S(b_Z)|\leq |\mm{V}(t_*,b)|/R(t_*,b_W)\leq M^2$.
    Decreasing the neighborhood $G'\rav W'\times Z'$ of $b$, we can assume $diam\, Z'<\gamma/M.$

 Fix $t'\in I.$  Consider $w'\in W',z',z''\in Z'$.  Set $y'\rav (w',z')\in G'.$
 Consider
 the following
 exit-time problem:
 \begin{eqnarray*}
  \left \{
 \begin{array}{rcl}
 \textrm{Minimize } \int_{t'}^T f_0(t,x,u)\, dt+\mm{V}(T;x(T))
 \\
 \textrm{subject to } \dot{x}=f(t,x,u),\quad u\in P, 
\\
x(t')=y',\quad T=\inf \{t>t'\,:\,x(t)\not\in\mm{W}\times \{z''\}\cup\{+\infty\} %
 \end{array} \right.
\end{eqnarray*}
 By the definitions of $Q_{\mm{W}}$ and $Q^{\mm{W}}$, for all $u\in\mm{U}$, there exists
 $T_u$ such that
 \begin{eqnarray*}
 t'+Q^{\mm{W}}(t',(w',z'),z'')\leq
 T_u&=&\inf \{t>t'\,:\,x_{y',t',u}(t)\not\in\mm{W}\times \{z''\}\}\\
 &\leq& t'+Q^{\mm{W}}(t',(w',z'),z''),
\end{eqnarray*}
 i.e.,
 \begin{eqnarray}\label{1047_1}
 0\leq T_u-t"\leq L||z'-z''||,\qquad x_{y',t',u}(T_u)\in\mm{W}\times \{z''\}.
\end{eqnarray}
  Since all $T_u$ are uniformly bounded  (by  $t'+L||z'-z''||$), the optimal value of this exit-time problem is $V(t',y');$
  in particular,
  $$V(t',y')\leq J(t',y';u,T_u)+V(T_{u},x_{y',t',u}(T_{u})).$$

  In addition, we can choose a control $u'\in\mm{U}$ with its motion $x'\rav x_{y',t',u'}$ such that
  $$|V(t',y')-J(t',y';u,T_u)-V(T_{u'},x'(T_{u'}))|\leq ||z'-z''||.$$
  By \rref{1047_1} we have $x'(T_{u'})\in \mm{W}\times \{z''\}.$ 
  Now, we can choose $\bar{w}\in\mm{W}$ such that
  $x'(T_{u'})=(\bar{w},z'').$ So,
 \begin{eqnarray}\label{1048_1}
  |R(t',w')S(y')-J(t',y';u,T_u)-R(T_{u'},\bar{w})S(z')|\leq ||z'-z''||.
\end{eqnarray}
  In addition, $T_{u'}-t'\leq L||z'-z''||\leq M\,diam\, Z'<\gamma.$ Then, by the choice of the compact $K$, we have
  $x'(t)\in K$ for all $t\in[t',T_{u'}].$ By the choice of  $M$, we obtain
 \begin{eqnarray}\label{1047_}
   ||\bar{w}-w'||&\leq& ||x_{y',t',u'}(T_{u'})-y'||\\
   &\leq& M(T_{u'}-t')\leqref{1047_1} M^2||z'-z''||\nonumber
\end{eqnarray}
and
  $|J(t',y';u,T_u)|\leq M(T_{u'}-t')$.
  Then,
   \begin{eqnarray}\label{1048_}
|R(t',w')S(z')-R(T_{u'},\bar{w})S(z'')|&\leqref{1048_1}& ||z'-z''||+M(T_{u'}-t')\\
&\leqref{1047_1}& 2M^2||z'-z''||\nonumber\\
&\leq& 2M\gamma\leq 2M^2.\nonumber
 \end{eqnarray}

   In particular, 
   in the case  $y'\rav y,z''\rav b_Z$, we can provide $(\bar{w},z'')\in K$, $\tau<\gamma$, and
  $$|R(t',w') S(z')|\leqref{1048_}
   |R(t'',\bar{w}) S(b_Z)|+2M^2\leq M |S(b_Z)|+2M^2\leq 3M^3$$
    by the choice of $M$.
   According  to $|R|>1/M$, we also have
   $$|S(z')|\leq 3M^4 \quad\forall (w',z')\in G'.$$

  Return to the general case.   Consider every  $w'\in W',z',z''\in Z'$  and, additionally, every $w''\in W',t\in I$. Then, we again obtain  $u'\in\mm{K}$, $T_{u'}\in[t',t'+\gamma]$ and $x_{y',t',u'}(T_{u'})=(\bar{w},z'')\in K$ satisfying \rref{1047_},\rref{1048_}. Now, by \rref{1048_},
  for all $t,t'\in I$, $w',w''\in W'$, and $z',z''\in Z',$
  $|\mm{V}(t',w',z')-\mm{V}(t,w'',z'')|$ does not exceed
\begin{eqnarray*}
   &\ & 2M^2||z'-z''||+               (R(T_{u'},\bar{w})-    R(t',w')+R(t',w')-R(t,w''))S(z'')\\
   &\leq& M^3||z'-z''||
   + 3M^5(T_{u'}-t'+||\bar{w}-w'||+||{w}''-w'||+|t'-t|)\\
   &\leqref{1047_}& M^3||z'-z''||+3M^5(1+M)(T_{u'}-t')+3M^5(||{w}''-w'||+|t'-t|)\\
   &\leqref{1047_1}& 7M^7||z'-z''||+3M^5(||{w}''-w'||+|t'-t|).
\end{eqnarray*}
 Thus,  $|\mm{V}(t',w',z')-\mm{V}(t,w'',z'')|\leq 7M^7(||z'-z''||+||w''-w'||+|t'-t|)$ holds
 for some neighborhood $I\times W'\times Z'$ of each $(t_*,b_W,b_Z)\in I\times W\times Z$.
  \qed

\section{Examples}

\begin{example}\label{ex22}
Let $x$ be the capital stock, $u$ the investment,  $\nu>0$ the depreciation rate, and $\mu$ the discount rate. Consider the following  problem:
\begin{eqnarray*}
 \textrm{Minimize } \int_{0}^{+\infty} e^{\mu t}g(x,u) dt,\\
 \textrm{subject to } \dot{x}=-\nu x+u,\ x(0)=x_*\quad u\in[0,U_{max}], \quad k>0.
\end{eqnarray*}
\end{example}

 Set $R(\theta)\rav e^{\mu \theta}$ for all $\theta\in\mm{R}.$ It is easy to prove that
$$J(\theta,y;u,T)=R(\theta)J(0,y;u',T-\theta)\qquad \forall  \theta\geq 0, y\in\mm{R}, T>\theta,$$ where $u'(t)\rav u(t+\theta)$ for all $t\geq 0.$
 Since
  $u\in\Omega$ iff $u'\in\Omega$ for each of
 $\Omega\in\{\Omega_{\mathcal{L}},\Omega_{\mathcal{R}},\mm{U}\},$
 we have $$\mm{V}(\theta,y)=R(\theta)\mm{V}(0,y)\qquad \forall y\in\mm{R}^m,\theta\geq 0,
 \mm{V}\in\{V^{\mathcal{R}},V^{inf},V^{\mathcal{L}},V^{all}\}$$ 
 if $\mm{V}$ is well-defined.

 Set $Z_<\rav]0,U_{max}/\nu[$.
 We obtain
 $0\in -\nu x + int\,[0,U_{max}]$ for all $x\in Z_<$.
 By Corollary~\ref{1112345},
 $\mm{V}$ is Lipschitz continuous (if finite) in $\mm{R}_{\geq 0}\times Z_<.$

 Set $Z_>\rav]U_{max}/\nu,\infty[$.
 We obtain
 $-\nu x + U_{max}<0$ for all $x\in Z_>$.
 By Corollary~\ref{11123456},
 $\mm{V}$ is Lipschitz continuous (if finite) in $\mm{R}_{\geq 0}\times Z_>.$


 So, the value functions are  Lipschitz continuous (if finite) in $\{(t,x)\in\mm{R}_{\geq 0}\times (\mm{R}_{\geq 0}\,|\,x\neq U_{max}/\nu\}$;  no assumptions on $g$ besides continuity are needed.

Consider a use case of Theoremae~\ref{11123},\ref{1112345}
 (compare with \cite{zelikin1994,Lev2015}):
\begin{example}\label{108}
\begin{eqnarray}
 \textrm{Minimize } \int_{0}^\infty g(y_1,y_2,u) dt\nonumber\\
 \textrm{subject to } \dot{y_1}=u,\ \dot{y_2}=y_1,\label{819}\\
  (y_1,y_2)(0)=b_*, \quad y_1,y_2,u\in P\subset \mm{R}.
  \nonumber
\end{eqnarray}
\end{example}
For a certain $k\in\mm{R}$, let us also require
$$ g(\nu y_1,\nu^2y_2,u)=\nu^k g(y_1,y_2,u)\qquad \forall \nu>0, (y_1,y_2,u)\in\mm{R}^{3}.$$

 Clearly, the functions from $\{V^{\mathcal{R}}, V^{inf}, V^{\mathcal{L}},V^{all}\}$ do not depend on $t$. More to come. First, note that $u\in\mm{U}$ iff all  maps $u'_\nu(t)\rav u(t/\nu)$ are within $\mm{U}$ for all $\nu>0.$
   Fix $(y_1,y_2)\in\mm{R}^{2},u\in\mm{U},\nu>0$. Set $x\rav x_{(y_1,y_2),0,u}$ and $x'\rav x_{(\nu y_1,\nu^2y_2),0,u'_\nu}$. For all $t\geq 0$, we have
 \begin{eqnarray*}
  x(t)=x'(t/\nu),\  g(x'(t/\nu),u(t/\nu))=\nu^k g(x(t),u(t)),\\
 J(0,(\nu y_1,\nu^2y_2);u'_\nu,T/\nu)=\nu^{k-1}J(0,(y_1,y_2);u,T).
 \end{eqnarray*}
 Now, for all finite $\mm{V}\in\{V^{\mathcal{R}}, V^{inf}, V^{\mathcal{L}},V^{all}\}$, we have
 $$\mm{V}(t,\nu y_1,\nu^2y_2)=\nu^{k-1}\mm{V}(0,y_1,y_2)\quad \forall \nu>0,t\geq 0,(y_1,y_2)\in\mm{R}^{2}.$$

Define $\mm{W}\rav\mm{R}_{>0}$, $Z\rav\mm{R}.$
 Let some $\mm{V}\in\{V^{\mathcal{R}}, V^{inf}, V^{\mathcal{L}},V^{all}\}$ be well-defined and finite.

 {\bf 
 The case $P=\mm{R}.$}

 Consider  $G_+=\{(y_1,y_2)\in\mm{R}^2\,:\,y_2>0\}$
 and the functions $R:\mm{R}_{\geq 0}\times\mm{W}\to \mm{R}_{>0}$ and $S:Z\to \mm{R}$
  defined as follows:
\begin{eqnarray}\label{1127}
 R(t,w)\rav w^{k-1}, S(z)\rav\mm{V}(0,z,1)\qquad   \forall z\in\mm{R},w>0,t\geq 0;
\end{eqnarray}
 then,
\begin{eqnarray}\label{1128}
 \mm{V}(t,y_1,y_2)=R(t,\sqrt{y_2})S(y_1/\sqrt{y_2})\quad \forall y_1\in\mm{R},y_2,t\geq 0.
\end{eqnarray}
  In the coordinates $(w\rav \sqrt{y_2},z\rav y_1/\sqrt{y_2})$, system \rref{819} has the  form
\begin{eqnarray}\label{1129}
\dot{w}=z/2,\quad \dot{z}=
  \frac{u-z^2/2}{w},\quad u\in \mm{R},\  (z,w)\in Z\times\mm{W}.
\end{eqnarray}
    Since for every point $(y_1,y_2)\in G_+$, in its sufficiently small neighborhood, the controls $u\rav\pm z^2(y_1,y_2)$ provide for~$z$ to increase/decrease with the speed at least  $1/4\,\sqrt{y_2}$.
    Set $L\rav 4\sqrt{y_2};$
     in view of Theorem~\ref{11123}, we find out that, in $\mm{R}_{\geq 0}\times G_+$
 (for $y_2>0$), the function $\mm{V}$ is Lipschitz continuous. The proof of the case $y_2<0$ is similar.

 Thus, in the case $P\rav\mm{R}$, the function $\mm{V}$ is Lipschitz continuous under
  $$\{(y_1,y_2)\in\mm{R}^2\,|\,y_2\neq 0\}.$$

 {\bf
 The subcase $P=[-a^2,a^2],$  $G^a_+\rav\{(y_1,y_2)\in\mm{R}^2\,:\,y_2>0,y^2_1<2y_2a^2\}.$}

 Consider the functions $R:\mm{R}_{\geq 0}\times\mm{W}\to \mm{R}_{>0},$ $S:Z\to \mm{R}$,
  defined by \rref{1127}. Then, $V$ satisfies \rref{1128},
    in the coordinates $(w\rav \sqrt{y_2},z\rav y_1/\sqrt{y_2})$, system \rref{819} has the  form
    \rref{1129}.
     Note that, for every point $(y_1,y_2)\in G^a_+$, we have $|y_1|<a\sqrt{y_2}$, i.e. $z(y_1,y_2)<a.$ Hence, in its sufficiently small neighborhood, the controls $u\rav\pm z^2(y_1,y_2)\in P$ provide for~$z$ to increase/decrease with the speed at least  $1/4\,\sqrt{y_2}$.
    Set $L\rav 4\sqrt{y_2};$
     in view of Theorem~\ref{11123}, we find out that, in $\mm{R}_{\geq 0}\times G^a_+$
 (for $y_1^2<a^2 y_2$), the function $\mm{V}$ is Lipschitz continuous. The proof of the case $-y_1^2>a^2 y_2$ is similar.

 {\bf 
  The subcase $P=[-a^2,a^2],$ $\widetilde{G}_+^a=\{(y_1,y_2)\in\mm{R}^2\,:\,y_2>0,y^2_1>2a^2y_2\}.$}


 Consider the functions $R:\mm{R}_{\geq 0}\times\mm{W}\to \mm{R}_{>0}$ and $S:Z\to \mm{R}$
  defined as follows:
 $$R(t,w)\rav w^{1-k}, S(z)\rav\mm{V}(0,1,z)  \qquad\forall z\in\mm{R},w>0,t\geq 0,$$
 then, $$\mm{V}(t,y_1,y_2)=R(t,y_1)S(y_2/y^2_1)  \qquad\forall y_1\in\mm{R},y_2,t\geq 0.$$
  In the coordinates $(w\rav y_1,z\rav y_2/y^2_1)$, system \rref{819} has the  form
$$
\dot{w}=u,\quad \dot{z}=\frac{1-2zu}{w},\  (z,w)\in Z\times\mm{W}.$$
     Note that, for every point $(y_1,y_2)\in \widetilde{G}_+^a$,
     we have $y^2_1>2a^2y_2$, i.e., $2a^2z(y_1,y_2)>1.$ Hence, in its sufficiently small neighborhood,
     for every admissible control, the coordinate $z$  strictly increases.
     Then,
     in view of Corollary~\ref{111234}, we find out that, in $\mm{R}_{\geq 0}\times \widetilde{G}_+^a$
 (for $y_1^2>2a^2 y_2$), the function $\mm{V}$ is Lipschitz continuous. The proof of the case $y_1^2<2a^2 y_2$ is similar.

Thus, in the case $P\rav[-a^2,a^2]$ the function $\mm{V}$ is Lipschitz continuous under
 $$\{(y_1,y_2)\in\mm{R}^2\,|\,y_2\neq 0, y_1^2\neq 2a^2\}.$$

\begin{example}\label{ex2}
 Let us refine Example \ref{ex22};
  the new example does not satisfy assumptions \cite[Hypothesis 3.1(iv)]{sagara}, \cite[Ch.3(A4)]{bc}, \cite[(C3)]{norv}, \cite[($\widetilde{A}6$)]{baumeister}.
\begin{eqnarray*}
 \textrm{Minimize } \int_{0}^{+\infty} \big[2u+|u|-x\big]\, dt\\
 \textrm{subject to } \dot{x}=2u-x,\ x(0)=b_*,\quad u\in P=[-1/2,1/2], \quad x\in\mm{R}.
\end{eqnarray*}
\end{example}
    Clearly, $J(0,b;u,T)=x_{b,0,u}(T)-b+||u||_{L_1([0,T],P)}$
    for all $b\in\mm{R},$ $u\in\mm{U}$, $T>0$. Since $x_{b,0,u}(T)>-2$ holds for large  $T$, the finiteness of the limit of $J(0,b;u,T)$ as $T\to\infty$ implies $u\in L_1(\mm{R}_{\geq 0},P)$, moreover,
  we have
  $$x_{b,0,u}(T)\to 0,\ b+J(0,b;u,T)\to||u||_{L_1(\mm{R}_{\geq 0},P)}\geq 0 \textrm{ as } T\to\infty.$$
  Then, we obtain $V^{inf}(0,b)=-b,$ i.e., $u^*\equiv 0$ is optimal in view of $V^{inf}.$
  Moreover, $V^{\mathcal{L}}\equiv V^{\mathcal{R}}\equiv V^{\inf}.$

  But, for $T>1$, we define a control $u_T\in\mm{U}$ as follows: $u_T(t)=-1/2$ for all $t\in [T-\ln 2, T]$, and
  $u_T(t)=0$ otherwise. Then, we have $x_{b,0,u_T}(t)=be^{-t}$  if $t\in[0,T-\ln 2]$, $x_{b,0,u_T}(t)=e^{-t}(b+e^T/2)-1$ if
  $t\in[T-\ln 2,T]$. Now, we obtain
  $$x_{b,0,u_T}(T)=be^{-T}-1/2, V^T(0,b)\leq J(0,b;u_T,T)=be^{-T}-1/2-b+ \frac{\ln 2}{2}.$$
   One easily proves that $V^T(0,b)=J(0,b;u_T,T).$ Passing to the limit, we obtain
    $V^{all}(0,b)=-b+ \frac{-1+\ln 2}{2}$ for all $b\in\mm{R}.$ Thus,
   $$V^{all}(b)<V^{\inf}(b)\quad\forall b\in\mm{R}.$$
   Since $V^{all}-V^{\inf}\equiv const$, \rref{optimal} guarantees that $u^*\equiv 0$ is optimal in view of $V^{inf}$ and $V^{all}$ at once. Moreover, $u^*\equiv 0$ is DH-optimal (in particular, agreeable)  and overtaking optimal
    but is not strongly optimal \cite{HaurieSethi}.

   Note that, in this example, the functions  $V^{all},V^{\inf}$ are well-defined and smooth.
   In \cite[Proposition 3.2]{motta2014}, for exit-time  control  problems,
 it was proved that $V^{all}=V^{\inf}$ if $V^{all}$ is well-defined and continuous. Moreover, in this case  $V^{all}=V^{\inf}$ is the unique nonnegative solution of the associated HJB equation. Conditions guaranteeing the continuity of $V^{all}$ were showed in \cite[Section~4]{motta2014}.

{\bf{Acknowledgements.}}
I would like to
express my gratitude to N.N.~Subbotina  and A.J.~Zaslavskii for valuable discussion in course of writing this article. Special thanks to  Ya.V.~Salii for the translation.

\end{document}